\documentclass[11pt,a4paper,twoside]{article}
     \usepackage{amsmath,amssymb,latexsym}
     
\newtheorem{prop}{Proposition}[section]
\newtheorem{thm}[prop]{Theorem}
\newtheorem{Def}[prop]{Definition}
\newtheorem{lemma}[prop]{Lemma}
\newtheorem{cor}{Corollary}[section]

\newcommand{\piu}[1]{+_{#1}}

\newcommand{\A}[1]{\mathrm{Aut}(#1)}

\newcommand{\gl}[2]{\mathrm{GL}(#1,#2)}
\newcommand{\pf}{\noindent{\bf Proof~}}
\newcommand{\ov}[1]{\overline{#1}}

\begin{document}
 \author{Carlo Casolo \\ \small
         Dipartimento di Matematica ``U.\,Dini",
         Universit\`a di Firenze\\ \small
         Viale Morgagni 67A, \ I-50134 Firenze, Italy\\ \small
         \texttt{casolo\,@\,math.unifi.it}\\[1ex]
         \and Orazio Puglisi\\ \small
         Dipartimento di Matematica ``U.\,Dini",
         Universit\`a di Firenze\\ \small
         Viale Morgagni 67A, \ I-50134 Firenze, Italy\\ \small
         \texttt{puglisi\,@\,math.unifi.it}}
         \title{\vspace*{-2cm}\bf
       Nil-automorphisms of groups with residual properties}

         \maketitle

%
%
%
%
\section{Introduction}

Following Plotkin \cite{plotkin} we say that the automorphism $x$ of the group $G$ is  a \emph{nil-automorphism} if, for every $g\in G$, there exists $n=n(g)$ such that $[g,_n x]=1$ (the commutator is taken in the holomorph of $G$ and, following the usual notation, the element
$[g,_n x]$ is defined inductively by $[g,_0,x]=g$ and, when $n>0$, $[g,_n x] =[[g,_{n-1} x],x]$).
If the integer $n$ can be chosen independently of $g$, then $x$ is said to be \emph{unipotent}.
Nil and unipotent automorphisms can be regarded as a natural extension of the  concept of Engel element, since a nil-automorphism $x$ is just a left Engel element in $G\langle x\rangle$. 
Another way to look at nil-automorphisms, is to consider them as a generalization of unipotent automorphisms of vector spaces.
For these reasons there are several natural questions that can be asked about nil-automorphisms, which are suggested by known facts
about Engel groups or unipotent linear groups.

The Engel condition is, however, rather weak and it is often necessary to add some additional hypothesis in order to investigate this property. Similar difficulties arise when we study groups of nil-automorphisms. In this paper we shall focus on groups satisfying some finiteness conditions.

 The first result we prove is a consequence of a classical theorem of Baer \cite{baer}. Recall that a  group $G$ satisfies Max if every subgroup of $G$ is finitely generated.

\bigskip
{\bf Theorem A}\emph{ If $G$ satisfies Max, and  $H$ is any group of nil-automorphisms of $G$, then   $[G,H]$ is contained in the Fitting radical of $G$ and $H$ stabilizes a finite subnormal series in $G$. In particular $H$ is  nilpotent.}

\bigskip

An interesting consequence of Theorem A is that groups of of nil-automorphisms of finite groups are nilpotent. On the other hand  every finite nilpotent group has a faithful representation as a group of nil-automorphisms of a suitable finite group, so that the following question seems  natural:

\bigskip

{\bf Question 1} \emph{ Is it true that every finite group of nil-automorphisms is nilpotent?}

\bigskip

Although the question is still far from a general solution, we provide a partial answer by restricting ourselves to a particular class of  groups. Namely the following theorem can be proved.

\bigskip

{\bf Theorem B} \emph{ Let $G$ be a locally graded group, and $H$ a finite group of nil-automorphisms of $G$. Then $H$ is nilpotent.}

\bigskip

The definition of locally graded groups will be recalled in the next section, but it is worth remarking that many interesting types of groups (e.g.  locally nilpotent groups, locally finite groups,  residually finite groups) are locally graded.

 A natural hypothesis, suggested by the ``linear" setting, is to require that for some  $n$, the identity $[g,_n x]=1$ holds for all $g\in G$ and $x\in H\leq \mathrm{Aut}(G)$. In this case we say that $H$ is a group of \emph{$n$-unipotent} automorphisms. It  is then reasonable to ask to what extent the linear case can be generalized.

\bigskip

{\bf Question 2} \emph{ Let $G$ be a  group and $H\leq \mathrm{Aut}(G)$. Assume that, for some $n$, the identity $[g,_n x]=1$ holds for all $g\in G$ and $x\in H$. Under which hypothesis on $G$  is  it true that $H$ is (locally) nilpotent?}

\bigskip

A well known result of Wilson \cite{wilson1} says that finitely generated residually finite $n$-Engel groups are nilpotent, 
suggesting that Question 2 may have  an affirmative answer for groups $G$ in this class. 
It is easy to see that, rather than looking at residually finite groups, it is  better to work with profinite groups, in 
order to take advantage of their well developed theory. We get the following theorem
\bigskip

{\bf Theorem C}\emph{ Let $G$ be a profinite group and $H$ an abstract finitely generated group of $n$-unipotent automorphisms. Assume that $G$ has a basis of open normal subgroups whose members are normalized by $H$. Then $F=cl_G([G,H])$ is pro-nilpotent and $H$ is nilpotent.
}

\bigskip

In particular theorem C applies when the profinite group $G$ is finitely generated because, in that case, $G$ has a basis whose elements are characteristic subgroups. Since any residually finite group can be embedded into its profinite completion, it is not hard, using theorem C, to derive a similar result for groups in this class.
 Another immediate corollary of  Theorem C is  the result of Wilson stated before.
 
Finally we wish to point to the reader's attention the recent work of Crosby and Traustason on Engel groups. In particular  the results contained in   \cite{traustason},  will be needed at some stage of our proofs.

\section{Groups with residual properties}

In order to investigate groups with residual properties, we need information about \linebreak nil-automorphisms of finite groups. This information will be obtained as a special case of a more general fact which, in turn, is a consequence of the following classical result of  Baer \cite{baer}.

\bigskip

\begin{thm}\label{baer} Let $G$ be a group satisfying Max. Then the set of left Engel elements of $G$ coincides with the Fitting radical of $G$. Consequently an Engel group with Max is a finitely generated nilpotent group.
\end{thm}

\bigskip

Theorem A can now be proved easily.

\bigskip

{\bf Theorem A}\emph{ If $G$ satisfies Max, and  $H$ is any group of nil-automorphisms of $G$, then   $[G,H]$ is contained in the Fitting radical of $G$ and $H$ stabilizes a finite subnormal series in $G$. In particular $H$ is  nilpotent.}

\pf Assume, first of all, that   $G$ is  abelian. Let $k$ be the rank of $G$ and set $T$ for the torsion subgroup of $G$. The subgroup $T$ is finite and $G/T$ is free abelian.  For every prime $p$, the group $A_p=G/G^p$ is elementary abelian and can be regarded as an $\mathbb F_p$-vector space of dimension at most $k$. Thus $H_p=H/C_H(A_p)$ is a unipotent group of automorhisms of  $A_p$, hence $[A_p,_k H_p]=1$. From this it follows that  $[G,_k H]\leq \bigcap_{p\in \mathbb P}G^p=T$.  On the other hand $T$ is finite and each element of $H$ acts as a nil-automorphism, so that there exists $m$ such that $[T,_m H]=1$. Thus $[G, _{k+m} H]=1$ and $H$ stabilizes a finite series in $G$. 
When $G$ is nilpotent, $H$  acts as a group of nil-automorphisms on each factor of the descending central series hence the above argument can be used to see that $H$ stabilizes a finite central series of $G$. 
Let now $G$ be any group satisfying Max, and choose $h\in H$.  It is clear that $h$ is a left Engel element of the group $G\langle h\rangle$ which, in turn, satisfies Max. By theorem \ref{baer} $h$ belongs to the Fitting radical of $G\langle h\rangle$. It follows that the subgroup  $[G,h]$ is contained in the Fitting radical of $G$. Therefore $[G,H]$ is contained in the Fitting radical of $G$ and property Max implies that $[G,H]$ is finitely generated, hence nilpotent. Now, by  the first part of the proof,  $H$ stabilizes a finite normal series in $[G,H]$, showing that $H$ stabilizes a finite subnormal series in $G$. By a well known result of Hall (see \cite{stability}), $H$ is nilpotent.  \hfill $\Box$

\bigskip 

As a corollary we obtain

\bigskip

\begin{cor}\label{finite} Let $G$ be a finite group and $H$ a group of nil-automorphisms of $G$. Then $H$ is nilpotent and $[G,H]\leq \mathrm{Fit}(G)$.
\end{cor}

\bigskip

We concentrate now on the main objects of our investigation, namely groups with residual properties. We start by recalling the definition of locally graded groups.

\bigskip

\begin{Def} A group $G$ is said to be \emph{locally graded} if every  non trivial finitely generated subgroup of $G$ has a proper normal subgroup of finite index.
\end{Def}

Of course any residually finite group is locally graded but the converse is not true. E.g. the groups $\mathrm{PSL}(2,\mathbb F)$, where $\mathbb F$ is an infinite locally finite field, are locally finite (hence locally graded) but, being simple, they do not have proper  subgroups of finite index.

\bigskip

Question 1 has a positive answer for locally graded groups.

\bigskip 

{\bf Theorem B} \emph{ Let $G$ be a locally graded group, and $H$ a finite group of nil-automorphisms of $G$. Then $H$ is nilpotent.}

\pf By way of contradiction we assume the claim false, so that there exist pairs $(G,H)$ where  $G$ is locally graded, $H\leq \mathrm{Aut}(G)$ is finite and consists of nil-automorphisms,  and $H$ is not nilpotent. Among them we select a pair $(G,H)$ in such a way that $H$ has minimal order. Therefore every proper subgroup of $H$ is nilpotent and the structure of $H$ can be easily described. There exist two different primes $p,q$, such that
$H=Q\langle x\rangle$ where $Q$ is an elementary abelian $q$-group, $x$ has order $p^m$ for some $m$,  and $Q$ is an irreducible $\langle x\rangle$-module.

In order to achive a contradiction, we shall prove that, for each $i\in \mathbb N$, if $[g,_i x]=1$ then $[g,Q]=1$.
When $i=0$ this is true, because $[g,_0 x]=g$. Assume we know the claim holds for all $j<i$, and choose $g\in G$ such that $[g,_i x]=1$. The element $a=[g,x]$ satysfies $[a,_{i-1} x]=1$, therefore $[a,Q]=1$. The subgroup $A=[g^{-1}, Q]$ is $H$-invariant. To see this we need only to show that $A$ is normalized by $x$. In fact for any $y\in Q$, we have
$$
[g^{-1},y]^x=[[x,g]g^{-1},y^x]=[a^{-1}g^{-1}, y^x]=[g^{-1}, y^x]\in A
$$

Assume $A\not=1$. The group $A$ is finitely generated, so that its finite residual $R$ is a proper subgroup. Since unipotent groups of automorphisms of finite groups are nilpotent, the subgroup $Q$ 
must act trivially on every finite image  given by $H$-invariant normal subgroups of $A$.
 Therefore $Q$ acts trivially on $A/R$. Consider the map
$\eta:Q\longrightarrow A/R$ defined by $(y)^\eta=[g^{-1},y]R$. If $u,v$ are in $Q$, then $[g^{-1}, v,u]\in R$, hence 
$$(uv)^\eta=[g^{-1}, uv]R=[g^{-1}, vu]R=[g^{-1}, u][g^{-1}, v]^uR=
[g^{-1}, u][g^{-1}, v][g^{-1}, v,u]R=(u)^\eta(v)^\eta
$$

The map $\eta $ is therefore a surjective homomorphism, hence $A/R$ is finite and $R$ is
 trivial. Thus $Q$ centralizes $A$. Moreover $A$ is isomorphic to a quotient of $Q$ and $x$
 acts on $A$ as a unipotent automorphism. Since $x$ has order $p^m$ the only possibility is 
 that $1=[A,x]=[g^{-1},Q,x]$. Pick $y\in Q$. Witt's identity gives
 $$
 1=[g^{-1},y^{-1},x]^y[y,x^{-1},g^{-1}]^x[x,g,y]^{g^{-1}}
 $$

On the other hand $[g^{-1},y^{-1},x]=1$ because $[g^{-1},Q,x]=1$, and $[x,g,y]=1$ because 
$Q$ centralizes $a=[g,x]$. Hence the above identity reduces to $[y,x^{-1},g^{-1}]=1$. 
This holds for every $y\in Q$, thus $[Q,x,g^{-1}]=1$. But $Q$ is an irreducible 
$\langle x\rangle$-module, so that $[Q,x]=Q$, proving that $Q$ centralizes $g^{-1}$, 
hence $g$. The inductive step is then complete.

Choose $g\in G$. There exists $n=n(g)$ such that  $[g,_nx]=1$.  Thus  $[g,Q]=1$ and  $Q$ turns out to  centralize $G$. Therefore $Q$ must then be trivial and this contradiction proves that the claim holds. \hfill $\Box$

\bigskip

We focus  now on Question 2. Since this question has a positive answer for finite groups, it seems  natural to study the same problem for  residually finite groups. 

To start off let us consider a finitely generated    residually finite group $G$. It is easily seen that there exists a residual system $\mathcal N$,  whose members are characteristic 
subgroups. The completion of $G$ with respect to $\mathcal N$ gives rise to a profinite group 
$\ov{G}$ on which $\mathrm{Aut}(G)$ acts as a group of continuous automorphisms. For this reason it is convenient to consider Question 2 in the class of  profinite groups and take advantage of the well-developed theory of profinite groups.

For each 
$N\in \mathcal N$ indicate by $N_0$ the closure of $N$ in $\ov{G}$. Each $N_0$ is normalized by 
$\mathrm{Aut}(G)$ so that there is  an induced action of  $\mathrm{Aut}(G)$ on $\ov{G}/N_0=GN_0/N_0\simeq G/N_0\cap G=G/N$.
If $H\leq \mathrm{Aut}(G)$ is a group of nil-automorphisms, then, for each $ N\in \mathcal N$,   its elements are nil-automorphims in their action on the finite group $\ov{G}/N_0 $. Therefore $[G,H]N_0/N_0$ is nilpotent, by corollary \ref{finite}. It readily follows that  the closure of $[G,H]$  in $\ov{G}$ is pro-nilpotent and it is isomorphic to the cartesian product of its Sylow $p$-subgroups. This argument suggests that  it might be useful, as a first step,  to consider Question 2 for pro-$p$-groups.

The first lemma we need is actually a fact about finite $p$-groups.

\bigskip

\begin{lemma}\label{nil on solv} Let $G$ be a solvable $p$-group of exponent $p^l$ and derived length $d$. If $H\leq\A{G}$ and $[g,_nh]=1$ for all $g\in G$ and $h\in H$, then $H$ has exponent bounded in terms of $p^l$, $n$ and $d$.
\end{lemma}

\pf The proof is by induction on the derived length of $G$. 

When   $G$  is abelian we can invoke Lemma 16 of \cite{casolo}. 
Assume now  the claim holds  for groups of derived 
length at most $d-1>1$, and let $G$ have derived length $d$. The group $K=H^{f(p^l,n,d-1)}$ 
acts trivially on both $G/G^{(d-1)}$ and $G^{(d-1)}$. Thus $[G,K,K]=1$ so that $[g,x^{p^l}]=[g,x]^{p^l}=1$ for all $g\in G, \, x\in K$. Therefore, 
setting $f(p^l,n,d)=f(p^l,n,d-1)p^l$, we get $[g,h^{f(p^l,n,d)}]=1$ for all 
$g\in G, \, h\in H$ proving that  $\mathrm{exp}(H)$ divides $ f(p^l,n,d)$, as  claimed.\hfill $\Box$

\bigskip

\begin{lemma}\label{nilaction} Let $H$ be a group of $n$-unipotent automorphisms of the abelian group $A$. If $[A,H]$ has finite exponent, then $H$ has finite exponent. 
\end{lemma}

\pf It is enough to prove the lemma when the exponent of $[A,H]$ is power of a prime $p$.
If $[A,H]$ has exponent $p$, define $k$ by $p^{k-1}\leq n<p^k$. Then, given $a\in A$ and 
$h\in H$, one has $1=[a,_nh]=a(h-1)^n=a(h-1)^{p^k}=a(h^{p^k}-1)$. Thus $h^{p^k}$ centralizes 
$A$ and $H$ has exponent dividing $p^k$. Arguing by induction assume that the  exponent of 
$[A,H]$ is  $p^{r}$. Fix $h\in H$ and consider its action on the group $A/[A,H]^{p^{r-1}}$. Since  $[A,H]/[A,H]^{p^{r-1}}$ has exponent $p^{r-1}$, there exists $m$ such that $H^m$ centralizes $A/[A,H]^{p^{r-1}}$. Thus $[A,H^m]\leq [A,H]^{p^{r-1}}$ which  has exponent $p$. Thus $H^m$ has exponent dividing $p^k$, and $H$ has exponent dividing $mp^k$.\hfill $\Box$

\bigskip

The next lemma shows that, when a group of automorphisms of a $p$-group $G$ is $n$-unipotent, then it is possible to construct several powerful subgroups of $G$.

\bigskip

\begin{lemma}\label{powerful1}Let $G$ be a
finite
$p$-group,
$H\leq \A{G}$ a
 group of $n$-unipotent automorphisms and $L$ any $H$-invariant subgroup of $G$.
Then   there exists $k$ such that
$R=[L,H^k]$ is powerful.
\end{lemma}

\pf We treat the case $p\not=2$  first. 
Let $k=f(p,n,2)$ as defined in  lemma \ref{nil on solv}, set $K=H^k$ and
$R=[L,K]$. The group $H$ acts on $R$ as a group of $n$-unipotent  automorphisms, hence
$[R,K]\leq R^{(2)}R^p$. Clearly $[R,K^g]\leq
R^{(2)}R^p$ for every element $g\in L$.  Therefore $R'=[R, [K,L]]\leq R^{(2)}R^p$.
Dedekind's modular law gives $R'=R'\cap R^{(2)}R^p=R^{(2)}(R^p\cap R')$ but, since
$R^{(2)}$ is contained in the Frattini subgroup of $R'$, we get $R'=R^p\cap R'$.
Therefore $R'\leq R^p$, proving the claim. When $p=2$ choose $k=f(4,n,2)$ and apply the same argument. \hfill $\Box$

\bigskip 

This lemma can be easily adapted to the case of pro-$p$-groups. Our main interest is the case of finitely generated pro-$p$-groups but it is better to consider a slightly more general setting, that we describe here below.

\begin{Def} We say that  $(G,H)$ is  a \emph{$(n,p)$-couple} if
\begin{enumerate}
\item $G$ is a pro-$p$-group and $H\leq\mathrm{Aut}(G)$ is finitely generated;
\item there exists $n$ such that each $h\in H$ acts as an $n$-unipotent automorphism on $G$;
\item $G$ has a basis of open subgroups whose members are normalized by $H$.
\end{enumerate}
\end{Def}

If $G$ is a finitely generated pro-$p$-group and $H$ is a finitely generated group of 
$n$-unipotent automorphisms of $G$, then $(G,H)$  is a $(n,p)$-couple. In fact, for any fixed natural number  $k$, $G$ contains only finitely many subgroups of index $k$, so that, given $N\trianglelefteq G$ open, the normal subgroup $N_H=\bigcap_{h\in H}N^h$ is $H$-invariant and of finite index. By a well known result of Serre, $N_H$ is open (see \cite{wilson}, theorem 4.3.5) hence, if $\mathcal B$ is a basis of open subgroups for $G$, then $\mathcal B_H=\{N_H\mid N\in \mathcal B\}$ is a basis of $H$-invariant open subgroups. 
On the other hand the set $\mathcal A_H=\{ C_H(G/N)\mid N\in \mathcal B_H\}$ is a basis of open subgroups for $H$, with respect to which $H$ is an Hausdorff topological group.

When considering a  topological group $U$, for every subgroup $A$,  we shall indicate by  $cl_U(A)$ its   topological closure  in $U$. As usual a closed subgroup $A$ is said to be \emph{finitely generated} when $A$ is the closure of a finitely generated discrete subgroup $A_0$.

\bigskip

\begin{lemma}\label{powerful2}Let $(G,H)$ be a  $(n,p)$-couple.  Then the following facts hold
 \begin{enumerate}
 \item There
exists
$k$ such that, for every  $H$-invariant subgroup $L$ of $G$,  the group  $R=cl_G([L,cl_H(H^k)])$ is powerful. 
\item If $d=d(H)$ is the minimal number of generators for $H$, there exists $r=r(p,n,d)$ such that, for every $x\in G$, the group  $L=cl_G(x^H)$  can be  generated by $r$ elements.
\end{enumerate}
\end{lemma}
  
\pf Let $k=f(p,n,2)$ or $f(4,n,2)$ when $p=2$,  as defined in  lemma \ref{powerful1}, and set
 $K=cl_H(H^k)$. 
Let $\mathcal B=\{N_\lambda\mid \lambda\in \Lambda\}$ be a basis of normal  $H$-invariant open subgroups of $G$ and choose $L$ an $H$-invariant subgroup of $G$.
 For each $\lambda\in \Lambda$ define  $L_\lambda=LN_\lambda/N_\lambda$ and 
 $R_\lambda=[L,K]N_\lambda/N_\lambda=[L_\lambda, K]$. 
 Since $K=\bigcap_{\lambda\in \Lambda}H^kC_H(G/N_\lambda)$, 
 we have $[L,K]\leq [L,H^k]N_\lambda$ so that $R_\lambda=[L,H^k]N_\lambda/N_\lambda=[L_\lambda,H^k]$. 
 Each $L_\lambda$ is
  a finite $p$-group and, by lemma \ref{powerful1},  $R_\lambda$ is powerful. Therefore  the group $R=cl_G([L,K])$, which is the inverse limit of the $R_\lambda$, is powerful.
 
The group $H/K$ is $d$-generated, residually finite  and of exponent bounded by $k$, which is a  function of $p$ and $n$. By  Zelmanov's solution of the restricted Burnside problem,  $H/K$ is finite of order bounded by a a suitable function  $r=r(p,n,d)$. 
Choose $x\in G$ and set  $L=cl_G(x^H)$. Arguing as in the first paragraph, we see that  $[L,K]\leq L^{(2)}L^p$ and, in particular 
$[L,K]\leq L'L^p\leq \Phi(L)$.
 Choose $X=\{h_i\mid i=1, \dots ,
r\}$ a left tranversal for
$K$ in $H$. Given any $y\in K$ and 
$h\in X$,  we have 
$$x^{hy}=x^h[x^h,y]\equiv
x^h\pmod{\Phi(L)}$$

Thus, for each $\lambda \in \Lambda$,   the set $X_\lambda=\{ x^hN_\lambda\mid h\in X\}$ 
generates  $L_\lambda=LN_\lambda/N_\lambda$, modulo its Frattini subgroup. Thus 
$\langle X_\lambda  \rangle=L_\lambda$. Since  each $L_\lambda$ can be generated by $r$
elements,  the same holds for their inverse limit $L$. 
\hfill
$\Box$

\bigskip

The next fact concerns the action of nil automorphisms on uniformly powerful group.

\bigskip

\begin{prop}\label{unipotent} Let $G$ be a finitely generated uniformly powerful pro-$p$-group and consider the usual additive structure defined on it. If $\alpha$ is a $r$-unipotent automorphism of $G$, then the automorphism induced by $\alpha$ on $(G,+)$ is $r$-unipotent.
\end{prop}

\pf Given $k\geq 1$ and $x\in G$, define 
$$
C(x,k)=x^{\alpha^k} 
x^{-{{k}\choose{1}}\alpha^{k-1}} 
\cdots 
x^{(-1)^{k-1}{{k}\choose{k-1}}\alpha}x^{(-1)^{k}}=
{\prod_{i=0}^k}x^{(-1)^i{{k}\choose{i}}\alpha^{k-i}}
$$
We stick to the notation of chapter 4 of \cite{an}. Fix $n\geq 1$ and any $x\in G$. The group $G_{n+1}/G_{2n+2}$ is abelian, whence 
$$[x^{p^n},_k \alpha]\equiv C(x^{p^n},k)\pmod{G_{2n+2}}
$$
and the relation
$$
(x^{p^n})^{(-1)^i{{k}\choose{i}}\alpha^{k-i}}=(x^{(-1)^i{{k}\choose{i}}\alpha^{k-i}})^{p^n}
$$
holds for every $i=0, 1, \dots , k$. Thus 

$$
(C(x^{p^n},k))^{p^{-n}}=x^{\alpha^k}
\piu{n} 
x^{-{{k}\choose{1}}\alpha^{k-1}}
\piu{n}
x^{{{k}\choose{2}}\alpha^{k-2}}
\piu{n} \dots 
\piu{n}
x^{(-1)^k}=a_n(x,k)
$$

Lemma 4.10 of \cite{an} shows that, for each $n$, 
$$a_n(x,k)\equiv b_n(x,k)=[x^{p^n},_k\alpha]^{p^{-n}}\pmod{G_n} $$
so that the  sequences $\{ a_n(x,k)\mid n\in \omega\}$ and $\{ b_n(x,k)\mid n\in \omega\}$ have the same limit. The limit of the first sequence is, by definition, the element 
$$\sum_{i=1}^kx^{(-1)^i{{k}\choose{i}}\alpha^{k-i}}=\sum_{i=1}^k (-1)^i{{k}\choose{i}}x^{\alpha^{k-i}}=[x,_k\alpha]$$
of  $(G,+)$. Since  $\alpha$ is $r$-unipotent, $[g,_r \alpha]=1$ for each 
$g\in G$. Thus each $b_n(x,r)=1$, so that $[x,_r\alpha]=0$ in $(G,+)$, showing that $\alpha$ 
acts unipotently on $(G,+)$.\hfill $\Box$

\bigskip

At this stage we can prove some facts  concerning unipotent automorphisms of nilpotent groups. These results will be used in the proof of the general case, but have some interest in their own. 
We begin by considering abelian groups.

\bigskip

\begin{prop}\label{action on abelian} Let $(A,H)$ be a $(n,p)$-couple and assume that $A$ is a torsion-free  abelian group. There exists $m=m(n)$ such that $[A,_m H]=1$.  In particular $H$ stabilizes a series of length at most $m$ in $A$ and it is nilpotent of class at most $m-1$.
\end{prop}

\pf We refer to the terminology defined in \cite{traustason}. Let $G=AH$ be the semidirect 
product of $A$ and $H$. Thus $A$ is a residually hypercentral $n$-Engel normal subgroup of $G$.
 Hence, by theorem 3 of \cite{traustason},  $A^{f(n)}\leq \zeta_{m(n)}(G)$ for suitable 
 functions $f,m$. Being $A$ torsion-free, it follows that $A\leq \zeta_{m(n)}(G)$, and $H$ 
 stabilizes a series of length at most $m(n)$ in $A$. \hfill $\Box$

\bigskip 

A straightforward consequence is the following.

\bigskip

\begin{cor}\label{uniform} Let $(G,H)$ be a $(n,p)$-couple and assume that $G$ is a  finitely generated uniformly powerful pro-$p$-group. Then $H$ is torsion-free and there exists $m=m(n)$ such that $H$ is nilpotent of class at most $m$.\end{cor}

\pf Consider $H$ as a subgroup of $\A{(G,+)}=\mathrm{GL}(r,\mathbb Z_p)$. By proposition \ref{action on abelian} $H$ is a unipotent subgroup of $ \mathrm{GL}(r,\mathbb Z_p)$ so that it is nilpotent and torsion-free. Clearly $((G,+),H)$ is a $(n,p)$-couple so that, again by \ref{action on abelian}, the nilpotency  class $H$ can be bounded in terms of  $n$. \hfill $\Box$

\bigskip

The following is an  easy fact that we shall need later. Its proof can be easily derived from the results about isolators collected in \cite{hall}.

\bigskip

\begin{lemma}\label{finite index}Let $V$ be a torsion-free nilpotent group,  $U$ a normal subgroup of finite index and $\alpha$ an automorphism of $V$.
\begin{enumerate}
\item The groups $U$ and $V$ have the same nilpotency class.
\item If $\alpha$ centralizes $U$, then $\alpha=1$.
\end{enumerate}
\end{lemma}

\bigskip

\begin{prop}\label{nilpotent1} Let $G$ be a finitely generated  nilpotent pro-$p$-group, and $H$ a $d$-generated group of $n$-unipotent  automorphisms of $G$. Then the following hold
\begin{enumerate}
\item  $H$ stabilizes a central series in $G$ and is therefore nilpotent;
\item  there exists $\kappa=\kappa(p,d,n)$ such that $\gamma_\kappa(H)$ is finite.
\end{enumerate}
 \end{prop}

\pf Under these hypotheses $G$ has an $H$-invariant basis of open subgroups, so that $(G,H)$ is a $(n,p)$-couple. 

Assume first  that $G$ is torsion-free. It is easily seen, arguing by induction on the nilpotency class of $G$, that $H$ is  torsion-free. The key observation is that a nil automorphism of a torsion-free abelian group, has finite order if and only if it is  the identity.
Let $K$ be the subgroup of $H$ described in lemma 
\ref{powerful2} and set $R=cl_G([G,K])$. The group $R$ is uniformly powerful because it is powerful and torsion-free (see theorem 4.8 in \cite{an}). The group $G$ acts on $R$ by conjugation, and this action gives rise to an action on $(R,+)$. Proposition \ref{unipotent} shows that each  element of $G$ acts unipotently on $(R,+)$, whence  $G/C_G(R)$ is isomorphic to a closed group of unitriangular matrices over $\mathbb Z_p$. Each non-trivial closed subgroup of  $\gl{r}{\mathbb Z_p}$ has rank bounded in terms of $r$ whence, by theorem 3.13 of \cite{an}, it possesses   a characteristic  powerful subgroup of finite index. Such subgroup is uniformly powerful, because the unipotent subgroups of $\gl{r}{\mathbb Z_p}$ are torsion-free.
Let  $\overline{M}=M/C_G(R)$ be such a subgroup. The group $H$ acts on both   $R$ and $M/C_G(R)$  and we can use corollary \ref{uniform} to see that   there exists $m=m(p,d,n)$ such that $\gamma_m(H)$ acts trivially on both $R$ and $M/C_G(R)$. Set $S=\gamma_m(H)\cap K$. Since $S\leq C_{GH}(R)$, then $[M,S]\leq [G,K]\cap C_G(R)=R\cap C_G(R)=\zeta(R)$ . Notice that, by proposition \ref{action on abelian},  $H$ stabilizes a finite series of length at most $m$ in $\zeta(R)$.
Fix $x\in M$ and consider the map $\theta_x: S\longrightarrow \zeta(R)$ defined by $(s)\theta_x=[s,x]$. It is readily seen that $\theta_x$ is an homomorphism. For any $k\in K, \, s\in S$ we have
$$
[s^k,x]=[s,x[x,k^{-1}]]^k=[s,x]^k
$$
because $[s,k^{-1}]$ belongs to $R$ and $S$ centralizes $R$. Thus $\theta_x$ is an homomorphisms
 of  $K$-modules, showing that $S/\ker(\theta_x)$ has a $K$-central series of length at most $m$. The subgroup $M$ has finite index in $G$, so that it is (topologically) finitely generated. Let $x_1, \dots , x_r$ be  a set of generators for $M$. Then 
 $$
 \bigcap_{i=1}^r \ker(\theta_{x_i})=\{ s\in S\mid [s,x_i]=1 \, \, \forall \, i=1, \dots , r\}=
 \{ s\in S\mid [s,x]=1 \forall \, x\in M\}=C_S(M)
 $$
 By lemma \ref{finite index} $C_S(M)=C_S(H)=1$.
For this reason   $S$ embeds, as a $K$-module,  into the direct product $\prod_{i=1}^r S/\ker(\theta_{x_i})$, so that $S$ has a finite $K$-central series of length at most $m$. This shows that $\gamma_{2m}(K)=1$ and,  by lemma \ref{finite index}, the same holds for $H$. Setting $\kappa(p,d,n)=2m$ gives the claim. 

For the general case set $T$ for the torsion subgroup of $G$ and put $\overline{G}=G/T$. If $G$ is abelian, then $T$ is pure in $G$, so that $T^p=T\cap G^p$. Thus $TG^p/G^p\simeq T/T\cap G^p=T/T^p$ and this shows that $T$ is finitely generated, hence finite. By induction on the nilpotency class of $G$, we see that $T$ is always finite and, in particular, closed in $G$.
Thus $\overline{G}$ is a torsion-free pro-group and , if $\kappa=\kappa(n,p,d)$ is the integer defined in the previous paragraph,  the group $\gamma_\kappa(H)$ acts trivially on $\overline{G}$. Any $h\in  \gamma_\kappa(H)$ is then completely determined, once  we know the $r$ elements $t_{h,i}=[x_i,h]$ so that   $\left| \gamma_\kappa(H)\right |\leq \left| T\right |^r$ and the claim is established.\hfill $\Box$

\bigskip

This result has a rather strong consequence.
\bigskip

\begin{prop}\label{nilpotent2}Let $(G,H)$ be a $(n,p)$-couple with $G$ a  torsion-free nilpotent pro-$p$-group. If   $H$ can be generated by  $d$ elements,  there exists $\kappa=\kappa(p,d,n)$ such that $\gamma_\kappa(H)=1$.
 \end{prop}

\pf Let $\kappa$ be the integer defined in proposition \ref{nilpotent1}. Given any $x\in G$, the group $L=cl_G(x^H)$ is, by lemma \ref{powerful2}, finitely generated. By proposition \ref{nilpotent1} we have $[L,\gamma_\kappa(H)]=1$, because $G$ is torsion-free so that a non-trivial  $n$-unipotent automorphism of $G$ can not have finite order. In particular, $\gamma_\kappa(H)\leq C_H(x)$. Being $x$ a generic element of $G$, this shows that $\gamma_\kappa(H)=1$.  \hfil $\Box$

\bigskip

When the prime $p$ is large enough, some useful bounds can be obtained.
\bigskip

\begin{prop}\label{bound} Let $G$ be a finite $p$-group of class $m$, and $H$ a group of $n$-unipotent automorphisms of $G$. There exists  integers $c=c(n,m), \, l= l(n)$ such that, if $p>l$, then $H$ is nilpotent of class at most $c$.
  \end{prop}
  
  \pf First of all assume $G$ abelian. Set $A=GH$ and let $k=k(n), \, l=l(n)$ be the integers  defined in lemma 3 of \cite{traustason}. For each $g\in G$ and $a\in A$, we have $[g,_n a]=1$ and  $[G,_r H]=1$ for some $r$. Lemma 3 of \cite{traustason} can be used to get that, for every $g\in G$ and $h_1, h_2, \dots , h_k\in H$, one has $[g, h_1, h_2, \dots, h_k]^l=1$. Since $p>l$ the element  
  $[g, h_1, h_2, \dots, h_k]$ must be trivial, showing that $[G, _kH]=1$. Thus   $[G,\gamma_k(H)]=1$ and the class of $H$ is bounded by $c(n,1)=k(n)-1$.
  When $G$ has class $m>1$ we consider the action of $H$ on the factors of the ascending (or descending) central series and apply the above argument. It is readily seen  that $[G, _{mk} H]=1$ so that $H$ has class at most $c(n,m)=mk(n)-1$. \hfill $\Box$

   \bigskip
   
   This proposition can be  immediately extended to pro-$p$-groups.
      
   \bigskip

   \begin{cor}\label{big prime} Let $(G,H)$ be a $(n,p)$-couple, with $G$ nilpotent of class $m$, and let $l(n), c(n,m)$ be the integers defined in proposition \ref{bound}. If $p>l(n)$, then $H$ is nilpotent of class at most $c(n,m)$   \end{cor}
   
 \pf If $\mathcal N$ is a basis of $H$-invariant open normal subgroups of $G$, then $G/N$ is a finite $p$-group of class at most $m$, for all $N\in \mathcal N$. The group $H$ acts on each $G/N$ as a group of $n$-unipotent automorphism so that, by proposition \ref{bound}, $\gamma_{c(n,m)+1}(H)$ centralizes all these quotients of $G$. Therefore $[G, \gamma_{c(n,m)+1}(H)]\leq \bigcap _{N\in \mathcal N}N=1$ and $H$ has class at most $c(n,m)$.\hfill $\Box$

   \bigskip

\bigskip

Before embarking in the proof of our main result, three  technical lemmata are needed. The first of them is due to Hartley and its  proof can be found in \cite{hartley}.

\bigskip

\begin{lemma}\label{carrier}Let $G$ be a group with an ascending series with abelian factors. If the members of the series are characteristic, then $G$ contains a characteristic subgroup $U$ such that $U$ is nilpotent of class at most two,  and  $C_G(U)=\zeta( U)$.
\end{lemma}

\bigskip

\begin{lemma}\label{double abelian} Let $G$ be a group and $H$ a  group of $n$-unipotent  automorphisms of $G$. Suppose  that $H$ has a normal abelian subgroup $A$ such that $H=A\langle h\rangle$,  and assume that $[G,A,A]=1$. Then $H$ is nilpotent.

\end{lemma}

\pf The subgroup $M=[G,A]$ is a normal abelian subgroup of $G$, and will be regarded as an $\langle h\rangle $-module. For each $k=1, 2, \dots, n$ let $X_k=\{g\in G\mid [g,_kh]=1\}$ and set $A_0=A, A_k=C_A(X_k)$ when $k=1, \dots , n$. The subgroups $A_k$ form a descending chain in $A$ and $A_n=1$.
  Choose $g\in X_1$ and consider the map 

  \begin{eqnarray}
  \nonumber
 \theta_g: A \longrightarrow & [G,A] \\
 \nonumber
 a\longmapsto & [a,g].
\end{eqnarray}

  An easy calculation shows that $\theta_g$ is an $\langle h\rangle$-module homomorphism, so that $A/\ker(\theta_g)$ is $\langle h\rangle$-isomorphic to a submodule of $[G,A]$. In particular $h$ acts as an $n$-unipotent automorphism on $A/\ker(\theta_g)$. Thus  
  $$
  \ov{A}=A_0/\bigcap_{g\in X_1}\ker(\theta_g)
  $$
  embeds, as an $\langle h\rangle$-module, into a cartesian product of copies of $M$, hence it 
  is acted upon by $h$ as an $n$-unipotent automorphim. Therefore the  element $h$  stabilizes a series of length at most $n$ in $A_0/A_1$. 
  
  Suppose we have already shown that $\langle h\rangle$ stabilizes a series of length at most $n$, in each factor $A_i/A_{i+1}$ for $i=1, \dots , k<n$.
   Let $S$ be the semigroup generated by $x=h^{-1}$, and $\mathbb ZS$ its semigroup algebra.     
  Pick $g\in X_{k+1}$ and define the map 
 \begin{eqnarray}
  \nonumber
 \theta_g: A_k \longrightarrow & M=[G,A] \\
  \nonumber
 a\longmapsto & [a,g].
\end{eqnarray}
  Using the fact that 
 $[g,h]\in X_k$, we see that, for every $a\in A_k$, 
 one has $(a^x)\theta_g=((a)\theta_g)^{x}$. Therefore $A_k/\ker(\theta_g)$ is isomorphic, as a  $\mathbb ZS$-module, to a submodule of $M$. In particular, being $(x-1)^n\in \mathbb ZS$, for every $u\in A_k/\ker(\theta_g)$, the identity $[u,_n x]=1$ holds. Then the same identity holds for the action on $A_k/A_{k+1}$, because this group embeds, as a $\mathbb ZS$-submodule, into a cartesian power of copies of $M$. Therefore $\langle h\rangle$ stabilizes a finite series of length at most $n$ in $A_k/A_{k+1}$. Hence $\langle h\rangle$ stabilizes a series of length at most $n^2$ in $A$, showing that $H$ is nilpotent. \hfill $\Box$

\bigskip 

\begin{lemma}\label{nilpotentbycyclic} Let $(G,H)$ be a $(n,p)$-couple and assume that $H=R\langle h\rangle$ with $R$ contained in the  Hirsch-Plotkin radical of $H$. Let  $N\trianglelefteq H$ be such that $H/N$ is nilpotent and $N\langle h\rangle$ is locally nilpotent. If $R$ is nilpotent then $H$ is nilpotent.
\end{lemma}

\pf The group $H$ is finitely generated and metanilpotent hence, by 15.5.1 of \cite{robinson}, the Hirsch-Plotkin radical of $H$ coincides with the Fitting subgroup and is nilpotent. The subgroup $N\langle h\rangle$ is subnormal in $H$ because $H/N$ is nilpotent. Thus $N\langle h\rangle$ is contained in $R$, the Hirsch-Plotkin radical of $H$. Hence $\langle h\rangle \leq R$ and from this it follows that $H=R$.\hfill $\Box$

\bigskip

\begin{prop}\label{solvable} Let $(G,H)$ be a $(n,p)$-couple,  and let $l(n), c(n,m)$ be the integers defined in proposition \ref{bound}. 
\begin{itemize}
\item[(1)] If $p>l(n)$, then $[G,\gamma_{c(n,2)+1}(H), \gamma_{c(n,2)+1}(H)]=1$ In particular  $\gamma_{c(n,2)+1}(H)$ is abelian.
\item[(2)] $H$ is nilpotent.
\end{itemize} 
\end{prop}

\pf Let $\mathcal N=\{N_i\mid i\in I\}$ be a base of $H$-invariant open subgroups of $G$
 and set $G_i=G/N_i$. Each $G_i$ is a finite $p$-group so that, by lemma \ref{carrier}, it has a 
 characteristic subgroup $U_i=M_i/N_i$ of nilpotency class at most two, such that 
 $C_{G_i}(U_i)=\zeta(U_i)$. The group $\ov{G}=\prod _{i\in I}G_i$ can be topologized by
  assigning the set  
  $\mathcal B=\{ \ov{N}_J=\prod_{j\in J}G_j\mid \left| I\setminus J\right| \textrm{ finite }\}$
   as a base of open neighborhoods for 
  the identity. In this way $\overline{G}$ is a compact pro-$p$-group, and $H$ acts on it as a 
  group of $n$-unipotent automorphisms. The subgroup $U=\prod_{i\in I}U_i$ is closed in 
 $\overline{G}$, is  acted upon continuously by $G$ and  $H$,  and $C_{\overline{G}}(U)=\zeta(U)$. 
Choose any cofinite subset $J$ of $I$ and define $\overline{U}_J=U\cap\overline{N}_J=\prod_{j\in J}U_j$. The subgroups 
$\overline{U}_J$ form a basis of open neighborhoods of the identity in $U$.

If $A=C_H(\overline{U}_J)$, we have that $[G,A]\leq C_G(\overline{U}_J)$. Thus  $[U_j,[G,A]N_j/N_j]=1$ for all $j\in J$, hence $[G,A]N_j/N_j\leq C_{G_j}(U_j)=\zeta(U_j)$ when $j\in J$. Hence 
$$[ [G,A]N_j/N_j,A]\leq [\zeta(U_j),A]\leq [U_j,A]=1$$
for all $j\in J$,
showing that $[G,A,A]\leq \cap_{j\in J}N_j=1$ which, in turn, proves that $A$ is abelian. In particular, if  $C=C_H(U)$,  then  $[G,C,C]=1$ and $C$ is abelian.
If $p>l(n)$  we can apply corollary \ref{big prime} to the action of $H/C_H(U)$ on $U$, showing that $\gamma_{c(n,2)+1}(H)\leq C_H(U)$ so that \emph{(1)} holds.

Let $X=\{h\in H\mid h^H \textrm{ is abelian }\}$.  Assume  that  
$U=\bigcup_{h\in H\setminus X}C_U(h)$.  Since $H$ consists of continuous automorphisms of $U$, 
each $C_U(h)$ is closed and, being $H\setminus X$ finite or countable, we may invoke Baire's theorem to
 prove that, for some $h\in H\setminus X$, the subgroup $C_U(h)$ contains an open subset. Therefore $C_U(h)$ must contain an open subgroup of the form    $\overline{U}_J$. 
The subgroup $\overline{U}_J$ is normalized by $H$, so that $\overline{U}_J\leq C_U(h^x)$ 
for all $x\in H$ showing that  $\overline{U}_J$ centralizes $h^H$. Thus $h^H\leq C_H(\overline{U}_J)$
 and, by the previous  paragraph, $h^H$ is abelian, a contradiction.  
 It is therefore possible to find $u\in U\setminus \bigcup_{h\in H\setminus X}C_U(h)$, so that $C_H(u)\leq X$. The group 
 $V=cl_U(u^H)$ is nilpotent of class at most $2$, it is    finitely generated ( see  lemma 
 \ref{powerful2}) and $C_H(V)$  is contained in $C_H(u)$. Since $C_H(V)$ is contained in $X$, for every element $c\in C_H(V)$ we have that $c^H$ is abelian and it follows that $C_H(V)$  is a $2$-Engel group. In particular $C_H(V)$ is nilpotent of class at most $3$. Moreover  proposition \ref{nilpotent1} shows that 
 $H/C_H(V)$ is nilpotent. Thus $H$ is metanilpotent and, in particular, it is solvable. 
 
 To complete the proof we argue by contradiction, selecting   $H$  a counterexample of minimal derived length. Under this assumption each finitely generated subgroup of $H'$ is nilpotent, so that $H'$ is contained in $R$, the Hirsch-Plotkin radical of $H$. However $H$ is metanilpotent hence, by   15.5.1 of \cite{robinson},  $R$ is nilpotent and coincides with the Fitting subgroup of $H$. Each counterexample of minimal derived length, is therefore nilpotent-by-abelian. 
 
 We may  choose  $H$ in such a way that $d(H/R)$, the minimal number of generators of $H/R$, 
 is minimal. If $H/R=\langle h_1R, \dots , h_lR\rangle  $ and $l\geq 2$, then the two subgroups
  $H_1=\langle h_1, R\rangle$ and $H_2=\langle h_2, \dots h_l, R\rangle$, are both normal and 
  nilpotent. To see this let $F$ be a finite subset of $H_i$ and consider $W=\langle F\rangle$.
   Clearly $\mathrm{Fit}(W)\geq W\cap R$, hence $W/\mathrm{Fit}(W)$ is an homomorphic image of 
   $ WR/R\leq H_i/R$, showing that $W/\mathrm{Fit}(W)$ has minimal number of generators strictly 
   smaller than $d(H/R)$. By assumption $W$ is nilpotent. The group $H_i$ is then locally 
   nilpotent and normal in $H$, thus it is contained in $R$, the Hirsch-Plotkin radical of $H$, 
   proving that $H=H_1H_2\leq R$, a contradiction.  Therefore $d(H/R)=1$ and it is possible to 
   write $H=R\langle h\rangle$ for a suitable $h\in H$. Moreover we may choose $H$ in such a way
    that the nilpotency class  of $R$ is minimal. Consider $S=C_H(cl_U(u^H))$ and the subgroup 
    $S\langle h\rangle$.   Since  $H/S$ is nilpotent, then $S\langle h\rangle$ is subnormal 
    in $H$ and, by lemma \ref{nilpotentbycyclic}, it can not be  locally nilpotent.  If, for every  $s\in S$, the subgroup $\langle s, h\rangle$ is
        nilpotent, it is easy to show that $S\langle h\rangle$ 
        is hypercentral, hence locally nilpotent. For this reason there must exist $a\in S$ such that $\langle a, h\rangle$ is not 
       nilpotent. Since $a\in S\subseteq X$, then $a^{\langle a,h\rangle}$ is abelian. 
 For this reason there is no  loss of generality in  
 assuming that $H=\langle a, h\rangle=A\langle h\rangle$ where $A=a^H$.
        Notice that $A\leq R$.

The group $H$ is then   metabelian, so that it satisfies Max-n, the maximal condition on normal subgroups. 

   If $W$ is a  pro-$p$-group, and 
  $\sigma:H\longrightarrow \mathrm{Aut}(W)$ is an homomorphism, we say that $\sigma$ is an 
  \emph{$n$-unipotent representation} if $(W,H^\sigma)$ is a $(n,p)$-couple.

  The set 
    $$
  \mathcal N=\{ N\mid N\leq H \textrm{ is the kernel of an } n\textrm{-unipotent representation and } H/N \textrm{ is not nilpotent}\}
  $$
  is not empty, so that it has a maximal element $M$. By looking at the representation related to $M$, we may assume that 
  \begin{itemize}
  \item if $1\not=N\trianglelefteq H$ and $H/N$ acts faithfully on a  pro-$p$-group $W$ in such a way that $(W,H/N)$ is a  $(n,p)$-couple,   then $H/N$ is nilpotent.
  \end{itemize}

Denote by $C$ the centralizer  of $U$ in $H$.

Two cases should be considered.

{\bf Case 1.} $C$ is not trivial.
  
   Thus $H/C$ is nilpotent and, as we have seen before,  $C$ is abelian and contained in $R$. Since $[G,C,C]=1$, we may apply lemma \ref{double abelian} to see that   the subgroup $C\langle h\rangle$ is nilpotent. Lemma \ref{nilpotentbycyclic} can now be used to see that $H$ is nilpotent, a contradiction.

 {\bf Case 2} $C=1$.

 The group   $H$ acts faithfully on $U$.
 If $U$ is torsion, then it must be of finite exponent, because it is nilpotent and compact. If this is the case, for each $u\in U$, the  subgroup $ cl(u^H)$  has, by lemma \ref{powerful2},  a uniformly bounded number of generators. By Zelmanov's solution of the restricted Burnside problem, $cl(u^H)$ has uniformly bounded order. In particular $cl(u^H)=u^H$ and, for some  $m$,   $H^m\leq C_H(u^H)$. Hence 
 $$
 H^m\leq \bigcap_{u\in U} C_H(u^H)=C_H(U)=1
 $$
 and $H$ is finite, hence nilpotent. Thus $U$ is not a torsion group. 
 
 When  $W$ is a nilpotent group, we shall indicate  by $T(W)$ its   torsion sugbroup.  Let $S=C_H(U/T(U))$ and, for each $s\in S$, consider the homomorphism
 
 \begin{eqnarray}
  \nonumber
 \psi_s: U \longrightarrow & U/\zeta(U) \\
  \nonumber
 u\longmapsto & [u,s]\zeta(U).
\end{eqnarray}

 This homomorphism is continuous and for this reason $[U,s]\zeta(U)/\zeta(U)$ is compact. On the other hand $[U,s]$ is contained in the torsion subgroup of $U$, hence it has finite exponent, say  $e_1$. For this reason $[cl(U^{e_1}),s]\leq \zeta(U)$. The same argument applies to the continuous homomorphism
 
 \begin{eqnarray}
  \nonumber
 \varphi_s: cl_U(U^{e_1}) \longrightarrow & \zeta(U) \\
  \nonumber
 u\longmapsto & [u,s]
\end{eqnarray}
proving that $cl_U(U^{e_1})^{e_2}\leq \ker(\varphi_s)$. It is then straightforward to see that, for each $s\in S$, there exists $e=e(s)$ such that $[cl_U(U^e),s]=1$. The subgroup $S$ can be generated, as a normal subgroup, by finitely many elements $s_1, \dots , s_l$. Choosing $e=\max\{e(s_i)\mid i=1, \dots, l\}$, we get $[cl_U(U^e),S]=1$. 
If $S\not=1$ then $P=C_H(cl_U(U^e))\geq S$ is not trivial  and therefore $H/P$ is nilpotent, 
by our choice of  $H$. Thus  $\gamma_a(H)\leq P$ for some $a$.
 The group $H$ acts on $W=U/cl_U(U^e)$, $(W,H/C_H(W))$ is a $(n,p)$-couple and $W$ is nilpotent of finite exponent.
 As we have seen before,  $H/C_H(W)$ is a finite $p$-group. There exists  $b$ such that
  $\gamma_b(H)\leq C_H(U/cl_U(U^e))$ and, if $c=\max\{a,b\}$, then 
  $D=\gamma_c(H)\leq C_H(U/cl_U(U^e))\cap C_H(cl_U(U^e))$. In particular $[U,D,D]=1$ and $D$ is abelian. By lemma \ref{double abelian} the subgroup $D\langle h\rangle$ is nilpotent and, by lemma \ref{nilpotentbycyclic}, $H$ would be nilpotent. For this reason the subgroup $S$ must be trivial.

 For every  $u\in U$ set $L(u)$ for the closure in $U$ of $u^H$, and $C(u)=C_H(L(u))$. Recall that each $L(u)$ can be generated by a uniformly bounded number of elements. If $C(u)=1$ for some $u$, then $H$ is nilpotent, by proposition \ref{nilpotent1}, a contradiction. Thus $C(u)$ is always non-trivial, and $H/C(u)$ is nilpotent,  because $L(u)$ is a  finitely generated pro-$p$-group, hence countably based.  
 By proposition \ref{nilpotent1}, there exists $\kappa$ such that $\gamma_\kappa(H)C(u)/C(u)$ is finite. Therefore $[L(u),\gamma_\kappa(H)]$ is contained in   $T(L(u))\leq T(U)$. It then follows  $\gamma_\kappa(H)\leq S=C_H(U/T(U))=1$, so that  $H$ is  nilpotent. This final contradiction proves that the claim holds so that  $H$ must be  nilpotent. \hfill $\Box$

  \bigskip

   We are now in a positon to prove  Theorem C.

   \bigskip

{\bf Theorem C}\emph{ Let $G$ be a profinite group and $H$ a finitely generated group of $n$-unipotent automorphisms. Assume that $G$ has a basis of open normal subgroups whose members are normalized by $H$. Then $F=cl_G([G,H])$ is pro-nilpotent and $H$ is nilpotent.
}

\pf  The fact that  $F$ is pro-nilpotent has been already noticed in the paragraph before lemma \ref{nil on solv}. Write $F=\prod_{p\in I}F_p$ where each $F_p$ is the $p$-Sylow subgroup of $F$. 
Assume first that $H$ acts faithfully on $F$. 
By theorem \ref{solvable} each $H/C_H(F_p)$ is nilpotent, say of class $c_p$   and, 
when $p>l(n)$, $[F_p,\gamma_{c(n,2)+1}(H), \gamma_{c(n,2)+1}(H)]=1$. Therefore, if 
$c=\max\{ c(n, 2), \, c_p \mid p\leq l(n)\}+1\}$, we get 
$[F,\gamma_c(H),\gamma_c(H)]=1$ showing that  $H$ is  abelian-by-nilpotent.

 If, by way of contradiction, we assume the claim false, the same argument used in the  proof 
 of theorem \ref{solvable} can be applied, to see that  there must exist a counterexample 
 $H$ that is nilpotent-by cyclic. Write 
$H=R\langle h\rangle$ with $R$ the Fitting subgroup of $H$.  Lemma \ref{double abelian}  
says  that $\gamma_c(H)\langle h\rangle$ is nilpotent, hence $H$ is   nilpotent, by lemma \ref{nilpotentbycyclic}, a contradiction. 

Therefore $H$ can not act faithfully on $F$ and 
 $C=C_H(F)$ is not trivial. We point out that the previous argument shows that $H/C$ is nilpotent. The group $C$ satysfies $[G,C,C]=1$, whence $H$ is abelian-by-nilpotent. Assuming the claim false we apply once again the argument of theorem \ref{solvable}, to show that $H$ may be assumed to be nilpotent-by-cyclic.  At this stage   the same considerations of the previous paragraph provide a contradiction, thus proving that $H$ must be nilpotent. \hfill $\Box$
 
\bigskip

As anticipated in the introduction,  our results imply a well known theorem due to Wilson \cite{wilson1}.

\bigskip

\begin{thm} Let $H$ be a finitely generated residually finite $n$-Engel group. Then $H$ is nilpotent
\end{thm}

\pf Let $G$ be the profinite completion of $H$. Since  $H/\zeta(H)$ is a finitely generated group of $n$-unipotent automorphisms of $G$, we apply theorem C to show that $H/\zeta(H)$ is nilpotent. Thus $H$ is nilpotent and the theorem is proved.\hfill $\Box$

\end{document}